\documentclass[a4]{amsart}

\usepackage{latexsym}
\usepackage{verbatim}

\newcommand{\Omegabar}{\overline{\Omega}}

\newcommand{\R}{{\mathbb R}}

\newcommand\Ha{\operatorname{Ha}^{(1)}}

\newcommand\WF{\operatorname{WF}'}

\newcommand\Cbill{C_{\operatorname{bill}}}

\newtheorem{theo}{{\sc Theorem}}[section]
\newtheorem{cor}[theo]{{\sc Corollary}}
\numberwithin{equation}{section}
\newcommand{\DL}{{\mathcal D}\ell(\lambda)}
\newcommand{\SL}{{\mathcal S}\ell(\lambda)}

\newtheorem{lem}[theo]{{\sc Lemma}}
\newtheorem{prop}[theo]{{\sc Proposition}}

\newenvironment{defin}{\medskip\noindent{\it Definition:\/} }{\medskip}

\newtheorem{maintheo}{{\sc Theorem}}

\newcommand\pa{\partial}

\newcommand\vol{\operatorname{vol}}

\newcommand{\RR}{{\mathbb R}}
\newcommand\Op{\operatorname{Op}}

\def\dbyd#1#2{\frac{ \partial #1}{\partial #2}}

\title[Billiards and boundary traces of eigenfunctions ]
{Billiards and boundary traces of eigenfunctions }

\author{Steve Zelditch}
\address{Department of Mathematics, Johns Hopkins University, Baltimore,
MD
21218, USA}
\email{ zelditch@math.jhu.edu}

\thanks{  \#DMS-DMS-0071358.}

\date{\today}

\begin{document}

\maketitle


\begin{abstract} This is a report on recent results with A.
Hassell on quantum ergodicity of boundary traces of eigenfunctions
on domains with ergodic billiards, and of work in progress with
Hassell and Sogge on norms of boundary traces. Related work by
Burq, Grieser and Smith-Sogge is also discussed.

\end{abstract}

\section{Introduction} This is a report on  recent work (and
work in progress) with A. Hassell and C. Sogge on boundary traces
 of eigenfunctions of the Laplacian on bounded domains and their
relations to the dynamics of the billiard map. Boundary traces of
eigenfunctions arise when one makes a reduction of the
interior/exterior  eigenvalue problem or wave equation of a domain
to its boundary. In the classical approach to boundary problems of
Neumann and Fredholm, the boundary reduction is made via layer
potentials and associated boundary integral operators. Our results
are based on  an analysis of these operators as semiclassical
Fourier integral operators.

We work on a Lipschitz  domain $\Omega \subset \R^n$ which  is
assumed to be compact and piecewise smooth with corners. We give a
precise definition below. The interior  eigenvalue problem  is:
\begin{equation} \left\{ \begin{array}{ll} - \Delta u_j = \lambda_j^2 u_j
\;\; \mbox{in}\;\; \Omega, & \;\; \langle u_j, u_k \rangle =
\delta_{jk} \\ \\
B u_j = 0 \;\; \mbox{on} \;\; \partial \Omega,
\end{array} \right. \end{equation} where
\begin{equation} B u_j =\left\{ \begin{array}{l}  u_j |_{\partial \Omega} \;
(\mbox{Dirichlet}),\;\;\mbox{resp.}\\ \\
\partial_{\nu} u_j + K u |_{\partial
\Omega},\;\; (\mbox{Neumann, Robin})\end{array} \right.
\end{equation}
Here, $\partial_{\nu}$ denotes the unit interior normal
derivative.

By the boundary trace $u_j^b$ of an interior eigenfunction $u_j$,
we mean the complementary operator to $B$, namely
\begin{equation} u_j^b (q) =  \partial_{\nu} u_j |_{\partial \Omega} \;\; (\mbox{Dirichlet}),\;\; \mbox{resp.}\;\; u_j^b (q) =
  u_j |_{\partial \Omega} \;\;\;\;\; (\mbox{Neumann, Robin}). \end{equation}

These boundary traces are a reduction to the boundary of the
eigenfuncions, and they become eigenfunctions of a reduction to
the boundary  of the resolvent, given by  the boundary integral
operator $F_{\lambda}$ on $\partial \Omega$ defined by
\begin{equation}
F_{\lambda} f(y) =\left\{ \begin{array}{ll} -2 \int_{\partial
\Omega} \frac{\pa}{\pa \nu_{y'}}
G_0(y,y',{\lambda}) f(y') dA(y'),  & Neumann\\ & \\
 2 \int_{\partial \Omega} \frac{\pa}{\pa \nu_{y}} G_0(y,y',{\lambda}) f(y') dA(y'), & Dirichlet \end{array}
\right. \label{Neumann-F}\end{equation} where $dA(y)$ is the
induced boundary area form,  and where
$$\begin{array}{lll}
G_0(z,z',{\lambda} ) & = & (2\pi)^{-n} {\lambda}^{n - 2} \int e^{i
{\lambda} (z-z') \cdot \xi} \frac1{ |\xi|^2 - 1 - i 0} \, d\xi
\\ & & \\ & = &

C \lambda^{n-2} (\lambda |z-z'|)^{-(n-2)/2} \Ha_{n/2-1}(\lambda |
z-z' |).
\end{array}$$ is the free outgoing Green's function $\RR^n$. Here, $\Ha_{n/2-1}$,
is the Hankel function.  Thus, the Dirichlet boundary integral
operator has kernel
$$\begin{array}{l} F_{\lambda}(y, y') = (2\pi)^{-n} {\lambda}^{n-1} \int e^{i {\lambda}(y-y') \cdot \xi}
\frac{\xi \cdot \nu_z}{ |\xi|^2 - 1 - i0 } \, d\xi.
\end{array}
$$
where $\nu_z$ is the unit inward pointing normal at $z$.  It is
elementary and classical (from Green's formula) that
\begin{equation} \label{FONE}
F_{{\lambda}_j} u_j^{b} =  \left\{ \begin{array}{ll} u_j^{b}, & Neumann \\ & \\
- u_j^b, & Dirichlet. \end{array} \right.
\end{equation}
One may think of $F_{\lambda}$ as
 a transfer operator or  effective Hamiltonian in the reduction to the boundary. More precisely,  as
 we will explain in Section \ref{BIO},   $F_h$ is the quantization of the billiard
 map $\beta$
 on the ball bundle $B^* \partial \Omega$, and is thus the reduction to the boundary of the quantum dynamics
of the wave group in the interior. Yet it is clearly a more
elementary operator than the Dirichlet and Neumann wave groups.

Our interest is in the relation between the dynamics of $\beta$ on
$ B^*\partial \Omega$ and the asymptotics of the boundary traces
$u_j^b$. In particular we study:

\begin{itemize}

\item Asymptotics of matrix elements
$\langle A(h_j) u_j^b, u_j^b \rangle$;

\item Asymptotics of $L^p$ norms: $||u_j^b||_{L^p}$.

\item Asymptotics of ratios $||u_j^b||_{L^p}/ ||u_j^b||_{L^2}$

\end{itemize}

The problems and results we present are to some extent  analogues
for boundary traces of results of \cite{SZ, ZZw} on eigenfunctions
on boundaryless manifolds or on interior eigenfunctions of
manifolds with boundary. In the next section, we state the main
results at this time.

Besides A. Hassell and C. Sogge, the author would like to  thank
M. Zworski for  helpful discussions about non-convex domains and
ghost orbits, and  D. Tataru for checking the statements of his
results in Section \ref{NORM} (they are stated differently in
\cite{T}).

\section{\label{RESULTS} Statement of results}

Let $\Omega$ be a bounded subdomain of $\RR^n$ with closure
$\Omegabar$.

\begin{defin}\label{cpm} We say that $\Omega \subset \RR^n$ is a piecewise smooth manifold  if the boundary $Y = \partial \Omega$ is strongly Lipschitz, and can be written as a disjoint union
$$
Y : = \partial \Omega = \;H_1 \cup \dots \cup H_m \cup \Sigma,
$$
where each $H_i$ is an open, relatively compact subset of a smooth
embedded hypersurface $S_i$, and $\Sigma$ is a closed set of
$(n-1)$-measure zero.

The sets $H_i$ are called boundary hypersurfaces of $\Omega$. We
call $\Sigma$ the singular set, and write $Y^o = Y \setminus
\Sigma$ for the regular part of the boundary.
\end{defin}

\subsection{\label{QEST} Statement of results on quantum ergodicity}

Our main results concern the quantum ergodicity of boundary
traces. To state the  results, we need some notation. In the
table, $\kappa$ denotes a $C^{\infty}$ function on $Y$ while  $k$
is the principal symbol of the operator $K \in \Psi^1(Y)$, and $d
\sigma$ is the natural symplectic volume measure on $B^*Y$. We
also define the function $\gamma(q)$ on $B^*Y$ by
\begin{equation}
\gamma(q) = \sqrt{1 - |\eta|^2}, \quad q = (y,\eta).
\label{a-defn}\end{equation}

\bigskip
\noindent \hskip 50pt {\large \begin{tabular}{|c|c|c|c|c|} \hline
\multicolumn{4}{|c|}{\bf Boundary Values} \\ \hline
 B & $Bu$  &   $u^{b}$ & $d\mu_B$   \\
\hline
 Dirichlet &  $u|_{Y}$ &  $\partial_{\nu} u |_{Y}$ & $\gamma(q) d\sigma $ \\
\hline
 Neumann &   $\partial_{\nu} u |_{Y}$  &  $u|_{Y}$  & $\gamma(q)^{-1} d\sigma$  \\
\hline
 Robin  &  $(\partial_{\nu} u - \kappa u) |_{Y}$ &  $u|_{Y}$ & $\gamma(q)^{-1} d\sigma$  \\ \hline
$\Psi^1$-Robin & $(\partial_{\nu} u - K u) |_{Y}$ & $u|_{Y}$ & $\displaystyle{\frac{\gamma(q)d\sigma }{ \gamma(q)^2 + k(q)^2}} $\\
\hline
 \end{tabular}}
\bigskip

\begin{maintheo}\label{main}
Let  $\Omega \subset \RR^n$ be a bounded piecewise smooth manifold
(see Definition~\ref{cpm}) with  ergodic billiard map. Let
$\{u_j^{b}\}$ be the boundary traces f the eigenfunctions
$\{u_j\}$ of $\Delta_B$ on $L^2(\Omega)$ in the sense of the table
above.
 Let $A_h$ be a semiclassical operator of order zero on $Y$. Then there is a subset $S$ of the positive integers, of density one,
such that
\begin{equation}\begin{gathered}
\lim_{j \to \infty, j \in S} \langle A_{h_j} u_j^b,
u_j^b \rangle = \omega_B(A), \quad B = \text{ Neumann, Robin or $\Psi^1$-Robin}, \\
\lim_{j \to \infty, j \in S} \lambda_j^{-2} \langle A_{h_j} u_j^b,
u_j^b \rangle = \omega_{B}(A), \quad B = \text{
Dirichlet},\end{gathered} \label{main-eqn}\end{equation} where
$h_j = \lambda_j^{-1}$ and $\omega_B(A) = \int_{B^* \Omega}
\sigma_A d\mu_B$.
\end{maintheo}

We note that the boundary traces are only normalized by
$||u_j||_{L^2(\Omega)} = 1$. We obtain new results on $L^2$ norms
by taking $A = I$. For the Neumann case, we have
\begin{equation} \label{L2D}
\lim_{j \to \infty, j \in S} \| u_j^b \|_{L^2(Y)}^2 = \frac{2
\vol(Y)}{\vol(\Omega)}, \end{equation} while for the Dirichlet
boundary condition they imply
\begin{equation} \label{L2N}
\lim_{j \to \infty, j \in S} \lambda_j^{-2} \| u_j^b \|_{L^2(Y)}^2
= \frac{2 \vol(Y)}{n\vol(\Omega)}. \end{equation}

To our knowledge, the only prior results on $L^2$ norms of
boundary traces are the upper and lower bounds of \cite{BLR, HT}
in the Dirichlet case, and the bounds implicit in \cite{T} in the
Neumann case. We review the latter in the last section.

\subsection{Statement of results on $L^{p}$ norms of boundary
traces}

We now consider $L^p$ norms of boundary traces for general
domains. We state a number of results of work in progress
\cite{HSZ}. For simplicity of exposition, we only state the
results for $L^{\infty}$ norms.

 First we consider estimates of
boundary traces and state general $L^{\infty}$ bounds analogous to
those which follow from the remainder estimate of
Avakumovic-Levitan for the spectral function of the Laplacian on a
boundaryless manifold (see \cite{SZ} for references). We denote by
\begin{equation} \label{SPROJ} E_{[a, b]}^b(x,x) = \sum_{j:
\lambda_j \leq \lambda \in [a,b] } |u_j^b(q)|^2
\end{equation} the boundary trace of  spectral projections kernel
(on the diagonal) for $\sqrt{\Delta}$ corresponding to the
interval $[a,b]$.

\begin{prop} We have:

$$\sum_{j: \lambda_j \leq \lambda}
|u_j^b(q)|^2 = \left\{\begin{array}{ll} C \lambda^{n + 2} +
O(\lambda^{n+1}) , & \mbox{Dirichlet}
\\ & \\\lambda^n + O(\lambda^{n-1}), & \mbox{Neumann}.
\end{array} \right.$$

Hence,

$$
||E_{[\lambda, \lambda + 1]}^b(x,x) ||_{L^{\infty}} =
\left\{\begin{array}{ll} O( \lambda^{n + 1}) , & \mbox{Dirichlet}
\\ & \\
 O( \lambda^{n - 1}), & \mbox{Neumann}.
\end{array} \right.$$

\end{prop}

Of course, the estimate holds for each term $||u_j^b(q)|^2
||_{L^{\infty}(\partial \Omega)}$ separately. We state it in the
above form because it is sharper than the statement for boundary
traces of eigenfunctions: for instance, the results implies the
same bound for boundary traces $\psi_j^b$ of quasimodes of order
$0$ satisfying $|| (\Delta - \lambda^2_j) \psi_j|| = O(1)$. We
also wish to  emphasize that our estimates on boundary traces of
eigenfunctions cannot be better than estimates for the above
spectral projections.

Even  for boundary traces of  eigefunctions, the above estimate on
$L^{\infty}$ norms is  sharp among all compact Riemannian domains
with boundary. Indeed, they are achieved on the northern
hemisphere of the standard $S^n$ by zonal spherical harmonics with
pole on the boundary: let $x_0$ be a point on the boundary of the
northern hemisphere (i.e. the equator), and let $Y^{N x_0}_0$ be
the $L^2$-normalized  zonal spherical harmonic of degree $N$ with
pole at $x_0$, i.e. the spherical harmonic which is invariant
under rotations fixing $x_0$. Then $\frac{1}{2} [Y^{N x_0}_0(x) +
Y^{N x_0}_0(x^*)]$ is a Neumann eigenfunction which is essentially
$L^2$-normalized and which has the same $L^{\infty}$ norm as $Y^{N
x_0}_0$. Here, $x^*$ is the reflection of $x$ through the
equation, i.e. $(x_1, \dots, x_n)^* = (x_1, \dots, x_{n-1}, -
x_n).$ Taking the normal derivative of the odd part of $Y^N_1$
under the reflection produces a similar sharp example in the
Dirichlet case. Eucliean examples will be discussed in the final
section.

However, our  main result on norms of boundary traces, as in the
interior case studied in \cite{SZ}, shows  that this result is
generically not sharp. In the following, the symbol $f(x) =
\Omega(g(x))$ is the negative of $f(x) = o(g(x))$, i.e. $f(x)$
grows at least as fast as $g(x)$ along a subsequence. In the
following, $E_{[\lambda, \lambda + o(1)]}$ is short for
$E_{[a,b]}$ where $a, b$ are the endpoints of a sequence of
shrinking intervals satisfying $a = \lambda, b = \lambda +
\epsilon(\lambda)$ where $0 < \epsilon(\lambda) = o(1)$.

\begin{maintheo} \label{RECURRENT} Suppose that
$$
||u_j^b||_{L^{\infty}}= \left\{\begin{array}{ll}
\Omega(\lambda_j^{\frac{n + 1}{2}}) , & \mbox{Dirichlet}
\\ & \\
 \Omega( \lambda_j^{\frac{n - 1}{2}}), & \mbox{Neumann}.
\end{array} \right.$$
or more generally that
$$
||E_{[\lambda, \lambda + o(1)]}(x,x)||_{L^{\infty}}=
\left\{\begin{array}{ll} \Omega(\lambda^{\frac{n + 1}{2}}) , &
\mbox{Dirichlet}
\\ & \\
 \Omega( \lambda^{\frac{n - 1}{2}}). & \mbox{Neumann}.
\end{array} \right.$$

Then there exists a recurrent point $q \in \partial \Omega$, i.e.
there exists $T > 0$ and a positive measure set of directions $\xi
\in S^*_q \Omega$ such that the billiard orbit $\exp_q (t \xi)$
returns to $q$ at time $T$.

\end{maintheo}

To put this another way, if the set of looping directions at $q$
has measure zero for all $q \in \partial \Omega$, then for any
$\epsilon > 0$, \begin{equation} \label{LIMSUP} \limsup_{\lambda
\to \infty} \lambda^{- \frac{n + 1}{2}} ||E_{[\lambda, \lambda +
\epsilon]}(q,q)||_{L^{\infty}} \to 0 \;\; \mbox{as}\;\; \epsilon
\to 0. \; . \end{equation}
 As pointed out in Section \ref{EXAMPLES},
 a recurrent point in the sense above cannot occur in a convex analytic Euclidean
domain. Hence,

\begin{cor} \label{CONVEX} If $\Omega \subset \R^n$ is a convex analytic domain
with Euclidean metric, then
$$
||u_j^b||_{L^{\infty}}= \left\{\begin{array}{ll}
o(\lambda_j^{\frac{n + 1}{2}}) , & \mbox{Dirichlet}
\\ & \\
 o( \lambda_j^{\frac{n - 1}{2}}), & \mbox{Neumann}.
\end{array} \right.$$

\end{cor}

The hypothesis in the Corollary can probably be weakened quite a
bit. For instance, it seems unlikely that real analytic Euclidean
domains of any kind can have recurrent points in our sense, and
the same may be true for smooth Euclidean domains (see also
  \cite{SV} for related conjectures).

  We note that although the result
on boundary traces appears similar to that of \cite{SZ} on
eigenfunctions on manifolds without boundary, there is a serious
difference in that the estimate in the boundaryless case is on the
ratio $||u_j||_{L^p}/||u_j||_{L^2}$ whereas in the case of
boundary traces, it is the interior eigenfunction and not the
boundary trace which is $L^2$ normalized. In the Dirichlet case,
 $||u_j^b||_{L^2}$ is known to be bounded above and below by a
 constant times $\lambda_j$ \cite{BLR, HT}, so one could restate the result in
 terms of the ratio $||u_j^b||_{L^{\infty}}/||u_j^b||_{L^2}$. In
 the Neumann case, $||u_j^b||_{L^2}$ can vary with the domain and
 metric.

\subsection{Remarks on methods and related results}

As mentioned above, our approach to boundary problems is to work
entirely on the boundary. In a recent article \cite{Burq}, Burq
takes the opposite route of `reducing' the study of ergodicity of
bounday traces  to the interior.  He extends the methods of
Gerard-Leichtnam \cite{GL} to reduce the proof of ergodicity of
boundary traces to the results of Zelditch-Zworski \cite{ZZw} on
ergodicity of interior eigenfunctions.  Both approaches have their
advantages, to which we devote a few remarks.

 Our motivation for working entirely on the boundary  comes from
several sources. First, it seems to us natural to prove results
about boundary traces by working with quantum dynamics on the
boundary, i.e. by making a semiclassical analysis of the operator
$F_h$. We also hope (and believe) that the semiclassical analysis
of $F_h$ and related operators has an independent interest and
other applications, for instance  to inverse spectral theory (see
\cite{Z1}).

A second motivation is that the
 use of  semiclassical  asymptotics of the
resolvent and of $F_h$ to relate spectrum and billiards is now the
standard approach in the physics literature on  eigenvalues of
bounded domains. The  method originated in the articles of
Balian-Bloch \cite{BB1, BB2} and has been developed in, for
instance, the articles \cite{AG, B, BS, BFSS, GP, THS, THS2, TV}.
A key reason for the wide use of the boundary reduction (as
explained to the author by A. Backer and R. Schubert) is that
boundary traces of eigenfunctions on smooth domains are much
easier to compute numerically than interior eigenfunctions.  See
the expository article \cite{B} by A. Backer
 (especially sections (3.3.1) - (3.3.7) and figure
$20$)  for a discussion of the numerical methods, for further
references,  and for pictures of boundary traces of
eigenfunctions. Despite the classical nature of layer potentials
and the induced boundary integral operators, and their common use
in physics and engineering, we are not aware of  prior
 mathematical studies of their semiclassical properties.

 The boundary integral method based on $F_h$ seems to us more elementary than the
 interior methods, but it does have the disadvantage that the
 phase of $F_h$, given by  the distance function between boundary
 points, is the generating function of both the interior and
 exterior billiard maps. Hence one has to deal with  the
 so-called ghost billiard orbits which have links outside the domain (see Section \ref{NONCVX}).
In the original version of
 \cite{HZ},  we  restricted to convex domains where ghost orbits
 do not arise. Using his interior method, Burq \cite{B} then proved
a more general
 ergodicity result for non-convex as well as convex domains. After
 that, in
 the second version of \cite{HZ}, we extended our methods and results to
general (possibly non-convex) domains. In the interim,
 discussions with M. Zworski were very helpful. The additional
 complication of ghost orbits of non-convex domains does not in
 the end turn out to be serious.

\section{\label{BIO} Boundary integral operators as semiclassical Fourier integral operators}

We begin by explaining the sense in which $F_h$ is the boundary
reduction or effective Hamiltonian for the resolvent. It is a
classical result of potential theory.

We recall that the single, resp. double layer potentials of a
domain $\Omega \subset \R^n$ are the operators
\begin{equation} \label{layers} \left\{ \begin{array}{l}
{\mathcal S} \ell({\lambda})f(x) =  \int_{\partial \Omega} G_0(x,
q, {\lambda}) f(q) dA(q), \\ \\  {\mathcal D} \ell({\lambda})f(x)
= \int_{\partial \Omega}\frac{\partial}{\partial \nu_y} G_0 x, q,
{\lambda}) f(q) dA(q),
\end{array} \right.
 \end{equation}
 where
$dA(q)$ is the surface area  measure on $\partial \Omega$, where
$\nu$ is the interior unit normal to $\Omega$, and  where
$\partial_{\nu} = \nu \cdot \nabla$. They induce the boundary
integral  operators (\ref{Neumann-F}).

Denote by $ R_{\Omega}^D(\lambda)$ (resp. $ R_{\Omega^c}^{N }
(\lambda) $) the resolvent of the interior Dirichlet Laplacian
(resp. the exterior Neumann Laplacian). Then the Fredholm-Neumann
reduction takes the form:
\begin{equation} \label{POT}\left\{  \begin{array}{l}  R_{\Omega}^D(\lambda) = 1_{\Omega} R_0(\lambda)1_{\Omega} + 1_{\Omega} {\mathcal
D}\ell(\lambda) (I + F_{\lambda})^{-1} {\mathcal S}(\lambda)^{tr}
1_{\Omega}, \\ \\
 R_{\Omega^c}^{N tr} (\lambda) = 1_{\Omega^c} R_0^{tr} (\lambda)1_{\Omega^c}
+ 1_{\Omega^c} \DL (I +  F_{\lambda})^{-1 } \SL^{tr} 1_{\Omega^c}
\end{array} \right.
\end{equation}
This formula has the form of a Grushin reduction of the Laplacians
to operators on the boundary in the sense of Sjostrand (see e.g.
\cite{DS}). Although we are not studying trace formula here, we
note that the combination of these two formula implies that (in
the sense of distributions),
$$\begin{array}{l}
Tr_{\R^2} [R_{ D}^{\Omega^c}(\lambda) \oplus R_{N }^{\Omega} -
R_{0 \rho}(\lambda)] =  \frac{d}{d \lambda} \log \det (I +
F_{\lambda}) ,.\end{array}$$ where $\det (I + F_{\lambda} )$ is
the Fredholm determinant. This formula is often used to determine
the eigenvalues and resonances of domains.

 It should be observed that the
Dirichlet  resolvent has the complexity of the Dirichlet wave
group, while $F_{\lambda}$ is no more complicated that the free
resolvent. It is $(I +F_{\lambda})^{-1}$ which has the complexity
of the Dirichlet resolvent or the Dirichlet-to-Neumann operator.
The relative simplicity of $F_{\lambda}$ is one of the attractions
of the boundary reduction.

\subsection{Quantized billiard map: Neumann } \normalsize
\bigskip

Since $F_{\lambda}$ is a reduction to the boundary of the
resolvent, it should be a kind of quantization of the billiard
map. We now explain in what sense this is true.

The {\it billiard map} $\beta : B^* Y^o \to \tilde T^* Y$ is
defined on the open ball bundle $B^* Y^o$ as follows: given $(y,
\eta) \in T^*Y$, with $|\eta| < 1$ we let $(y, \zeta) \in S^*
\Omega$ be the unique inward-pointing unit covector at $y$ which
projects to $(y, \eta)$ under the map $T^*_{\pa \Omega} \Omegabar
\to T^* Y$. Then we follow the geodesic (straight line) determined
by $(y, \zeta)$ to the first place it intersects the boundary
again; let $y' \in Y$ denote this first intersection.
 If $y' \in
\Sigma$ then we define $\beta(y, \eta) = y'$. Otherwise,  let
$\eta'$ be the projection of $\zeta$ to $T^*_{y'}Y$. Then we
define
$$
\beta(y, \eta) = (y', \eta').
$$
The map $\beta_{-} : B^* Y^o \to \tilde T^* Y$ is defined
similarly, following the backward billiard trajectory (that is,
the straight line with initial condition $(y, 2 (\zeta \cdot
\nu_{y}) \nu_{y} - \zeta )$).

The  graph
\begin{equation}
\Cbill = \operatorname{graph} \beta \equiv \{ (\beta(q), q) \mid q
\in \mathcal{R}^1 \}. \label{billiardCR}\end{equation} of $\beta$
is a smooth Lagrangian submanifold of $B^* Y^o \times B^* Y^o$. In
 a neighbourhood of $(y_0, \eta_0, y'_0, \eta'_0)$ it is given by
\begin{equation} \Cbill =
\{ (y, -\nabla_{y} d(y, y'), y',  \nabla_{y'} d(y,y') ),
\label{local-d}\end{equation} where $d(y, y')$ is the Euclidean
distance function.
 For
strictly convex $\Omega$ it is given globally by \eqref{local-d},
for $y, y' \in Y^o$,  but this is not true in general. This causes
extra difficulties for nonconvex domains, namely it introduces
spurious billiard orbits (known as ghost orbits in the physics
literature)  which do not remain entirely in the domain, and which
have to be shown to be irrelevant.

We see directly from (\ref{Neumann-F}) that $F_h$ is an
oscillatory integral operator with phase equal to the generating
function $d(y, y')$ of $\Cbill$. This is especially clear in
dimension $3$ when
$$
G(x, x', \lambda) = \frac1{4\pi} \frac{e^{i\lambda |x - x'|}}{|x -
x'|}.
$$
However, it is also clear that $F_h(y, y')$ has the homogeneous
singularity on the diagonal of a pseudodifferential operator of
order $-1$, and this is how it is usually described. At least in
the convex case, one may extend $\beta$ to the boundary $|\eta| =
1$ of $B^*Y$ by fixing $S^*Y$ pointwise.  Then the diagonal of
$B^*Y \times B^*Y$ and $\Cbill $ intersect in the diagonal of
$S^*Y \times S^*Y$, so one may view the wave front set of $F_h$ as
being the union of   two intersecting Lagrangeans.  $F_h(y, y')$
thus has similarities to the the oscillatory integral operator
kernels associated to two intersecting Lagrangeans of
 Melrose-Uhlmann. However, the intersection occurs at the boundary
 of both Lagrangeans and is thus outside the scope of the class of
 Melrose-Uhlmann operators, presumably explaining the unusual
 composition law described in \cite{Z1}.
Due to the explicit nature of our problem, we  carry out the
analysis by hand without making use of general operator theories.
The main point to observe is that the Fourier integral part is of
order $0$ while the pseudodifferential part is of order $-1$, so
the Fourier integral part dominates.

\begin{prop} Assume that $ \partial \Omega$ is smooth.  Let $U$ be any
 neighbourhood of $\Delta_{S^*  \partial \Omega}$. Then there is a decomposition of $F_{\lambda}$ as
$$
F_{\lambda} = F_{1,{\lambda}} + F_{2,{\lambda}} + F_{3,{\lambda}},
$$
where $F_1$ is a Fourier Integral operator of order zero
associated with the canonical relation $\Cbill = graph (\beta)$,
$F_2$ is a pseudodifferential operator of order $-1$ and $F_3$ has
operator wavefront set contained in $U$.
\end{prop}

This implies an  Egorov type result for the operator
$F_{\lambda}$:

\begin{prop} \label{Egorov} Let $A_h = \Op(a_h)$ be a zeroth order operator whose symbol
 $a(y,\eta,0)$ at $h=0$ is supported away from  $|\eta| = 1$. Put $h = \lambda^{-1}$ and let $\gamma$ be
 as above, and
 let $\beta$ denote the billiard ball map on $B^*Y$. Then
$$
F_{\lambda}^* A_h  F_{\lambda} = \tilde A_h + S_h,
$$
where $\tilde A_h$ is a zeroth order pseudodifferential operator
and $\| S_h \|_{L^2 \to L^2} \leq C h$. The symbol of $\tilde A_h$
is
\begin{equation}
 \tilde a = \begin{cases}
\gamma(q) [ \gamma^{-1} (\beta(q)) a(\beta(q))] , \quad q \in B^*Y \\
0 , \phantom{\gamma(q)^{-1} \gamma(\beta(q)) b(\beta(q)) , \ } q
\notin B^*Y. \end{cases} \label{Egorov-formula}\end{equation}
\end{prop}

This is a rigorous version of the statement that $F_h$ quantizes
the billiard ball map.

\section{Boundary quantum ergodicity}

Quantum ergodicity is concerned with  quantizations of classically
ergodic Hamiltonian systems. It is essentially a {\it convexity}
result relating  the {\it time average}
$$\langle A \rangle := \lim_{T \to \infty} \frac{1}{2T}
\int_{-T}^T U_t A U_t^* dt,$$ and the (constant) {\it space
average},
$$\omega(A) \; I$$
of an observable.  Here, we use the notation of the interior
problem where $U_t$ is the wave group. Also,  $\omega(A)$ is an
invariant state on the algebra of observables, which arises from
the local Weyl law:
\begin{equation}
\lim_{\lambda \to \infty} \frac{1}{N(\lambda)} \sum_{\lambda_j
\leq \lambda}  \langle A  u_j, u_j \rangle = \omega(A).
\end{equation}

The system is  {\it quantum ergodic} if
\begin{equation} \label{QE} \langle A \rangle = \; \omega(A)\; I\; + \; K,
\;\;\mbox{with}\;\;\; \frac{1}{N(\lambda)} || \Pi_{\lambda} K
\Pi_{\lambda}||_{HS} \to 0, \end{equation} where $||\cdot||_{HS}$
is the Hilbert-Schmidt norm and where $\Pi_{\lambda}$ is the
spectral projection onto the span of eigenfunctions of
$\sqrt{\Delta}$ of eigenvalue $\leq \lambda.$ In terms of the
$L^2$-normalized eigenfunctions $u_j$, QE says:
$$\frac{1}{N(\lambda)} \sum_{j: \lambda_j \leq \lambda} |\langle A
u_j, u_j \rangle - \omega(A)|^2 \to 0,$$ for any observable $A$.
Following the work of Schnirelman, Colin de Verdiere, Zelditch in
the boundaryless case, it was proved by Gerard-Leichtnman
\cite{GL} ($C^{1,1}$ domains with Dirichlet boundary conditions)
and Zelditch-Zworski \cite{ZZw} (general case with corners) that
domains with ergodic billiards are quantum ergodic.

\subsection{Proof of Theorem \ref{main}}

The main steps in the proof are:
\begin{itemize}

\item The local Weyl law;

\item Analysis of the classical limit state and the $L^2$ ergodic theorem on the classical level.

\item A convexity inequality to convert quantum ergodicity to the
classical  $L^2$  ergodic theorem.

\end{itemize}

\subsubsection{Local Weyl law}

It  states:

\begin{prop} \label{LWL}
 Let $A_h$ be a semiclassical operator of order
 zero on $\partial \Omega$.   Then for any of the above boundary conditions $B$,
 we have:
 \begin{equation}
\lim_{\lambda \to \infty} \frac{1}{N(\lambda)} \sum_{\lambda_j
\leq \lambda}  \langle A_{h_j} u_j^{\flat}, u_j^{\flat} \rangle =
\omega_B(A) , \quad h_j = \lambda_j^{-1},
\label{Weyl-limit}\end{equation}
 where $\omega_B$ is the
state defined in the table in Section \ref{QE}.
\end{prop}

  In the case of
multiplication operators, the Weyl law was first proved by S.
Ozawa \cite{Ozawa}.

\subsubsection{Mean ergodic theorem for the classical limit
states}

It is obvious that the states
\begin{equation} \rho_j^{\flat}(A) :=  \langle A_{h_j} u_j^{\flat}, u_j^{\flat} \rangle
\end{equation}
are invariant for $F_{{\lambda}_j}$:
\begin{equation} \rho_j^{\flat}(F_{{\lambda}_j}^* A F_{{\lambda}_j}) = \rho_j^{\flat}(A)
\end{equation}
Here, as above, $h_j = \lambda_j^{-1}$.

Since  any average or limit of averages of these states  will be
invariant, the  local Weyl law implies:
\begin{cor} \label{INV} The states $\omega_B$ are invariant under $F_{\lambda}$: $\omega_B(F_{\lambda}^* A F_{\lambda}) =
\omega_B(A)$. \end{cor} The formulae for the limit states in the
table in Section \ref{QEST}  are found by explicit calculation of
traces, e.g. of the dual heat or wave traces. The limit measures
may be understood as follows:

\begin{itemize}

\item The Neumann limit measure $\frac{d \sigma}{\gamma(q)}$ is the
projection to $B^* \Omega$ of the interior  Liouville measure
restricted to the set $S^*_{\in}
\partial \Omega $ of inward pointing unit vectors to $\Omega$
along $\partial \Omega$ under the map $\pi(x, \xi) = (x, \xi^T)$
taking a (co)-vector to its tangential component.

\item The Dirichlet limit measure $\gamma(q) d \sigma$
arises because boundary values of eigenfunctions involve the
normal derivative in both factors before restricting to the
boundary.
    The symbol of $h \pa_\nu$, restricted to the spherical normal bundle, and then projected to $B^* \partial \Omega$
    is equal to $\gamma$, so we expect to get the square of this factor in the Dirichlet
    case compared to the Neumann case.

    \end{itemize}

We now consider the mean ergodic theorem for these classical limit
states. Let us first recall the result for the standard `Koopman'
operator associated to the billiard map:
$$T: L^2(B^*\partial \Omega, d \sigma)
\to L^2(B^*
\partial \Omega, d\sigma),\;\;T f (\zeta) =  f(\beta (\zeta)).$$
 From the invariance it follows that $T$
is a unitary operator.  When $\beta$ is ergodic, the unique
invariant $L^2$-normalized eigenfunction is a constant $c$, and
one has the mean ergodic theorem
\begin{equation} \label{UNITARY} \lim_{N \to \infty} ||\frac{1}{N}
\sum_{n = 1}^N T_F^n(f) - \langle f, c \rangle || \to 0.
\end{equation}

But the transformation provided by Proposition \ref{Egorov} (
Egorovs theorem) is  defined by:
$$T f (\zeta) = \frac{\gamma (\zeta)}{\gamma(\beta(\zeta))} f(\beta(\zeta)).$$ This $T$ is {\it not}  unitary
on $L^2(B^*\partial \Omega, d \sigma)$. The invariance
$\omega_B(A) = \omega_B(F_h^* A F_h)$ implies that the associated
measure  $d \mu_B$ is invariant under
$$T^* f (\zeta) = \frac{\gamma (\beta(\zeta))}{\gamma(\zeta)} f(\beta(\zeta)).$$
Simple calculations show:

\begin{itemize}

\item (i) The unique positive $T^*$-invariant density is given by $\gamma^{-1}
d\mu.$  The unique positive $T$-invariant density is given by
$\gamma  d\mu.$

\item (ii)  $T$ is unitary relative to the   inner product $\langle \langle, \rangle \rangle$ on $B^* \partial \Omega$
defined by the measure $d\nu = \gamma^{-2} d \mu$.

\end{itemize}

 When $\beta$ is ergodic, the   orthogonal projection
$P$ onto the $T$-invariant $L^2$-eigenvectors has the form $$ P(f)
= \langle \langle f, \gamma \rangle \rangle = [\int_{B^*
\partial \Omega} f \gamma \gamma^{-2} d \mu_B] \gamma
$$
Thus, $P(\sigma_A) = [\int_{B^*
\partial \Omega} \sigma_A \gamma^{-1} d \mu_B] \gamma =  \omega_B(A) \gamma$.

The mean ergodic theorem thus says:
$$\frac{1}{N} \sum_{k = 0}^{N-1} T^k \sigma_A \to \omega_B(A)
\gamma. $$

\subsubsection{Convexity}

To show that
$$
\langle A u_j^{\flat}, u_j^{\flat} \rangle \to \omega_B(A),
$$
along a density one subsequence of integers $j$ is essentially to
show that
\begin{equation}\label{DESIRED}
\limsup_{\lambda \to \infty} \frac1{N(\lambda)} \sum_{\lambda_j <
\lambda}  \Big| \langle (A -  \omega_B(A)) u_j^{\flat},
u_j^{\flat} \rangle \Big|^2 = 0. \end{equation}

Due to the novel form of the local Weyl law and the Egorov
theorem, we first prove an auxiliary result of this kind for the
quantization of the invariant symbol
\begin{equation} \begin{array}{l}  \sigma_R(q) \sim c \gamma(q) = c (1 - |\eta|^2)^{1/2},\\ \\
\end{array} \label{c}\end{equation}
with $c$ a normalizing constant.

\begin{lem} \label{R} For all $\epsilon > 0$, there exists a pseudodifferential operator
$R$ of the form (\ref{c}) such that
$$
\limsup_{\lambda \to \infty} \frac1{N(\lambda)} \sum_{\lambda_j <
\lambda}  \Big| \langle (A -  \omega_B(A) \; R) u_j^{\flat},
u_j^{\flat} \rangle \Big|^2 < \epsilon. \label{want-small} $$
\end{lem}

We prove this intermediate step by the usual time-average +
convexity:
$$\begin{array}{l}  \frac1{N(\lambda)} \sum_{\lambda_j <
\lambda}  \Big| \langle (A -  \omega_B(A) \; R) u_j^{\flat},
u_j^{\flat} \rangle \Big|^2 \\ \\=   \sum_{\lambda_j < \lambda}
\Big| \langle (\langle A_h \rangle_N -  \omega_B(A) \; R )
u_j^{\flat}, u_j^{\flat} \rangle \Big|^2,\\ \\\leq C
\sum_{\lambda_j < \lambda}   \langle |\langle A_h \rangle_N
 -
\omega_B (A)\; R |^2 u_j^b, u_j^b \rangle,\end{array}$$ where $
\langle A_h \rangle_N = \frac1{N} \sum_{k=1}^N ((F_h^k)^* A_h
F_h^k. $ By the local Weyl law,  the limit equals $\omega_B((A_N -
\omega_B(A) R )^2) $, or
\begin{equation}
= \int_{B^* \partial \Omega} \Big( \langle \sigma_A \rangle_N -
\omega_B(A) \sigma(R) \Big)^2 d \mu_B.
\label{BR3-intro}\end{equation}

By Proposition \ref{Egorov},  the symbol $\sigma_{\langle
A\rangle_N}$ of $\langle A\rangle_N$ is   $\frac{1}{N} \sum_{k =
0}^N T_F^k(\sigma_A)$. By the mean ergodic theorem, this converges
to
$$\begin{array}{l} P(\sigma_A) = c \gamma(q)  \times \int \sigma_A \gamma^{-1} \,
d\sigma  = c \, \omega_B(A) \;  \gamma(q) \end{array}.$$  This is
approximately equal to the symbol of $\omega_B(A) R$. Thus,
$$  \int_{B^* \partial \Omega} \Big( \sigma_{A_N} - \omega_B(A) \sigma(R) \Big)^2 d
\mu_B$$becomes small as $N \to \infty$. Thus, $\frac1{N(\lambda)}
\sum_{\lambda_j < \lambda}  \Big| \langle (A - \omega_B(A) \; R)
u_j^{\flat}, u_j^{\flat} \rangle \Big|^2$ is   arbitarily small as
$N \to \infty$, proving Lemma \ref{R}.

Finally, we need to go from  $$\frac1{N(\lambda)} \sum_{\lambda_j
< \lambda}  \Big| \langle (A - \omega_B(A) \; R) u_j^{\flat},
u_j^{\flat} \rangle \Big|^2$$  to
 $$\frac1{N(\lambda)} \sum_{\lambda_j < \lambda}  \Big|
\langle (A - \omega_B(A) \; I) u_j^{\flat}, u_j^{\flat} \rangle
\Big|^2.$$

It suffices to show that
 $\frac1{N(\lambda)} \sum_{\lambda_j < \lambda}  \Big|
\langle ( I  -  \; R) u_j^{\flat}, u_j^{\flat} \rangle \Big|^2$ is
arbitrarily small. But this is the case $A = I$ above, since
$\omega_B(I) = 1.$

\subsubsection{\label{NONCVX} Non-convex domains}

We now adapt the argument to nonconvex domains. The problem with
the argument above  is that the canonical relation of
$F_{\lambda}$  is larger than the billiard relation, since it
relates points on the boundary which are connected by a straight
line even if the line passes outside the domain (ghost orbits).
However, one can modify $F_{\lambda}$ so that it leaves the
boundary traces invariant and so that its wave front relation has
an arbitrarily small measure of ghost orbits. That this should be
possible was emphasized to us by M. Zworski.

As is easy to see, the property $F_{\lambda_j} u_j^b = u_j^b$
(which follows from Green's formula) is valid for any choice of
Green's function on a neighborhood of $\Omega$.  So we  modify the
metric outside the domain, while keeping it Euclidean inside. Let
$b$ be a smooth, compactly supported nonnegative function on
$\RR^n$ which vanishes on $\Omegabar$. Consider the metric
$$
g_s = (1 + sb) g_{\operatorname{Euclidean}} \ \text{ on } \RR^n.
$$
For sufficiently small $s$, no geodesics of $g_s$ starting at a
point in $\Omegabar$ have conjugate points in $\Omegabar$. Let
$G_s(z, z, {\lambda}) = (\Delta_s - ({\lambda} + i0)^2)^{-1}(z,
z')$ denote the kernel of the outgoing resolvent of the Laplacian
on $\RR^n$ with respect to the metric $g_s$. It has the parametrix
$$
G_t(z, z', {\lambda}) = {\lambda}^{n-2} e^{i {\lambda}
\operatorname{dist}_s(z, z')} a(z, z', {\lambda}), \quad z \neq
z'.
$$
We then define:
$$
F^s_{\lambda}(y, y') = 2 \dbyd{}{\nu_y} G_s(y, y', {\lambda}),
\quad y, y' \in Y,
$$
and average over $s$ to obtain the operator
\begin{equation}
\tilde F_{\lambda} = \int_0^1 \chi(s) F^s_{\lambda} \, ds.
\label{average-op}\end{equation} Here $\chi$ is smooth, supported
in $ (0, \delta)$, and nonnegative with integral $1$. This
averaging removes all but an arbitrarily small measure of ghost
orbits or spurioius points from the $\WF(F_{\lambda})$, i.e.
points$(y, \eta, y', \eta')$, with $\eta = d_y |y - y'|$, $\eta' =
-d_{y'} |y - y'|$ such that the line $\overline{yy'}$ leaves
$\Omegabar$. We still have
$$
\tilde F_{\lambda_j} u_j^b = u_j^b.
$$
and in addtion the averaged operator has, as in the convex case,
a decomposition
$$
\tilde F_{\lambda} = \tilde F_{1,{\lambda}} + \tilde
F_{2,{\lambda}} + \tilde F_{3,{\lambda}},
$$
where $\tilde F_1$ is a Fourier Integral operator of order zero
associated with the canonical relation $\Cbill$, $\tilde F_2$ is a
pseudodifferential operator of order $-1$ and $\tilde F_3$ has
operator wavefront set contained in $U$. Moreover,   the symbol of
$\tilde F_{1,{\lambda}}$ is the same as for $F_{\lambda}$. Thus,
the proof  quantum ergodicity in the non-convex case now runs as
in the convex case, with this averaged operator in place of
$F_{\lambda}$.

\section{\label{NORM} Norm estimates of boundary traces}

We now turn to the results  on $L^p$ norms of boundary traces.
They are derived from an analysis of the singularities of the
boundary trace
\begin{equation} \label{BTW} E_B^b(t, q, q) =  \sum_{j = 1}^{\infty} e^{i t \lambda_j}
|u_j^b(q)|^2 \end{equation} of the  kernel of the wave operator
$\cos t \sqrt{\Delta_B}$ on  $\Omega$ with boundary condition $B$.
The nice feature of the boundary trace of $E_B(t, x, x)$ is that
the singularity at $t = 0$ becomes uniformly isolated from other
singularities, while the interior kernel $E(t, x, x)$ has
singularities at $t = 2 d(x)$ arbitrarily close to $t = 0$.  We
note that $E_B^b(t, q, q')$ is a spectral transform of the
Dirichlet-to-Neumann kernel.

We recall that the boundary traces $u_j^b$ are normalized by
$||u_j||_{L^2(\Omega)}$, so that it is of interest even to obtain
estimates on  $||u_j^b||_{L^2(\partial \Omega)}$. In
(\ref{L2D})-(\ref{L2N}) we obtained asymptotics of a density one
sequence of boundary traces in the ergodic case. In general, the
results are:

\begin{itemize}

\item Under a non-trapping assumption, $C_1 \lambda \leq ||u_j^b||_{L^2(\partial \Omega)} /||u_j||_{L^2(\Omega)} \leq C_2 \lambda$ in the
Dirichlet case \cite{BLR, HT};

\item For any smooth domain and metric, $||u_j^b||_{L^2(\partial \Omega)} /||u_j||_{L^2(\Omega)} = O(\lambda_j^{1/3}) $
in the Neumann case
  (Tataru, \cite{T}, Theorem 3);

 \item For a Sinai billiard on a torus (the exterior of a convex
 smooth domain in the torus),  $||u_j^b||_{L^2(\partial \Omega} /||u_j||_{L^2(\Omega)}
 = O(\lambda_j^{1/6}) $ in the Neumann case (\cite{T}, Theorem 5).

\end{itemize}

We introduce some notation. Given $(q, \eta) \in B^*_q \partial
\Omega$, we let $\xi(q, \eta) \in S^*_{in, q} \Omega$ be the
inward pointing unit vector to $\Omega$ obtained from $(q, \eta)$
by adding a multiple of the unit normal. We further denote by
$\Phi^t$ the broken bicharacteristic flow of the wave group, i.e.
the flow which carries singularities of the solution of the wave
equation. Finally, we let $v^T$ denote the tangential projection
of $v \in S^*_q \Omega$ at $q \in \partial \Omega.$ Finally, we
denote by $\gamma^B_q$ the boundary trace in the $q$ variable
corresponding to the boundary condition $B$.

\subsection{Singularity of boundary trace of wave kernel at $t =
0$}

Isolation of  the singularity at $t = 0$ of $E_B^b(t, q, q)$
follows from simple wave front set considerations. We first note
that
$$WF (\gamma_q^B \gamma_{q'}^B E(t, q, q')) \subset \{(t, \tau, q,
\eta, q', \eta'): [\Phi^t(q, \xi(q, \eta))]^T = (q', \eta'), \;
\tau = - |\xi| \}, $$ which  follows from propagation of
singularities for the wave equation. Composition  with the
boundary trace just pulls back (i.e. restricts) this wave front
relation to the boundary.

It follows that  $$WF (\gamma_q^B \gamma_{q'}^B E(t, q, q))
\subset \{(t, \tau, q, \eta, q, \eta'): [\Phi^t(q, \xi(q,
\eta))]^T = (q, \eta'),\;\; \tau = - |\xi(q, \eta)|\}. $$

Thus, the  singularities of the boundary  trace $\gamma_q^B
\gamma_{q'}^B E(t, q, q)$ at $q \in \partial \Omega$  to broken
bicharacteristic
 loops based at $q$ in
$\overline{\Omega}$. Let us first consider loops in convex
domains. They are either:

\begin{itemize}

\item closed geodesics on $Y = \partial \Omega$;

\item $m$-link transversal reflecting rays with $m$ vertices in $Y$, one of which is at $q$ and
which satisfy Snell's law of equal angles except possibly at $q$.

\end{itemize}

In the non-convex case, the description is more complicated since
there exist additional gliding rays and since boundary geodesics
need not carry singularities. For instance, for Euclidean plane
domains in  dimension $2$, bicharacteristics enter and exit the
boundary at inflection points and glide only over the convex part
of the boundary.

By definition of normal singularity, there exists $r \in R_+$ such
that $t^r E^b(t, q,q) \in C^{\infty}$ near $t = 0$.  Thus,  there
exist coefficients $a_j(q)$ such that
\begin{equation} E^b(t, q, q) \sim t^{- r} \sum_{j = 0}^{\infty}
a_j(q) t^j \end{equation} The next step is to calculate the
coefficients of the singularity.

Using a method of Ivrii, we first prove:

\begin{prop} \label{CONORMAL} The singularity of $ E^b(t,x,x)$
at $t = 0$ is a classical conormal singularity. \end{prop}

Granted the proposition, we can calculate the coefficients (and
the order) of the singularity using Hadamard type variational
formulae for eigenvalues:
$$ \delta \lambda_j = \int_Y \rho |u_j^b(y)|^2
dA(y),$$ where $\delta$ denotes a variation of either the boundary
or the boundary conditions and $\rho$ represents the tangent
vector to the variation. Thus, we can determine the integrals
$\int_{\partial \Omega} \rho(q) a_j(q) dA(q)$ and hence the
coefficients $a_j(q)$  asymptotically by considering the variation
of the wave trace \begin{equation}\label{it}  \frac{1}{it} \delta
\; \sum_{j} e^{i t \lambda_j} =\; \sum_{j} (\delta \lambda_j )
e^{i t \lambda_j} = \sum_{j} e^{i t \lambda_j} \int_Y \rho
|u_j^b(y)|^2 dA(y). \end{equation}

\subsubsection{Dirichlet}

In this case, we  vary the boundary in the normal direction with
variation vector field $\rho \nu$. The wave trace formula has the
form:
$$ \sum_{j} e^{i t \lambda_j} = C_n Vol_n(\Omega) (t + i 0)^{-n} +
C_{n-1}  Vol_{n-1}(\partial \Omega) (t + i 0)^{-n + 1} + a_1(t + i
0)^{-n + 2 } +\cdots  $$ where $\cdots$ represents lower order
terms (cf. \cite{I}, Theorem 2.1 or  \cite{M}, Corollary (3.6)).
It follows that \begin{equation} \label{DSING}     \sum_{j} e^{i t
\lambda_j} \int_Y \rho |u_j^b(y)|^2 dA(y) = C_n \delta
Vol_n(\Omega) (t + i 0)^{-n-1 } + C_{n-1} \delta
Vol_{n-1}(\partial \Omega) (t + i 0)^{-n } + \cdots \end{equation}
Here, we divided by the coefficient $it $ in (\ref{it}). The
variation of $Vol(\Omega)$ is non-zero, so we get a $(t + i 0)^{-n
- 1}$ term on the right.

\subsubsection{Neumann boundary conditions}

In this case, we vary the boundary conditions, and therefore
consider more general boundary conditions of the form
$$\partial_{\nu} \phi(q)  + \kappa(q) \phi(q) = 0, \; q \in \partial \Omega.$$
We denote by $\delta $ the first variation relative to a change in
the boundary condition from $\kappa \to \kappa + \epsilon \rho.$
By the Hadamard variational formula (studied by S. Ozawa for
general boundary conditions),
 the jth Neumann eigenvalue $\mu_j$ has
the variation
$$\delta \mu_j = \int_{\partial \Omega} |u_j^b(q)|^2 \rho(q)
d\sigma(q). $$

We now consider the first variation of the Neumann wave trace
under variations of pseudodifferential boundary conditions. We
obtain:
\begin{equation} \delta Tr E(t) = i t  \sum_j \delta \mu_j e^{i t \mu_j}. \end{equation}
The  $\rho$ term  only influences the third term of the
singularity trace expansion at $t = 0$ for Neumann-Robin boundary
conditions (see \cite{GM} for the calculation) Hence the variation
has the form
\begin{equation}  \delta Tr E(t) =   C_n (\int_{S^*\partial \Omega}
 \rho d \mu)  (t + i 0)^{-n + 2 } + \cdots, \end{equation}
from which we conclude
\begin{equation}  \label{NSING} \sum_j e^{i t \mu_j} [\int_{\partial \Omega} \rho |u_j^b|^2 dA]    =
C_n  (\int_{S^*\partial \Omega} \rho d \mu)  (t + i 0)^{-n + 1 } +
\cdots \end{equation}

\subsection{Spectral asymptotics of boundary traces}

We apply standard Tauberian theorems to obtain spectral
asymptotics of boundary traces from the singularity of the
boundary trace of the wave kernel.

\begin{prop} \label{LWL} We have:
$$\sum_{j: \lambda_j \leq \lambda}
|u_j^b(q)|^2  = \left\{\begin{array}{ll}  C \lambda^{n + 2} +
O(\lambda^{n+1}), & \mbox{Dirichlet}
\\ & \\
 C \lambda^n + O(\lambda^{n-1}), & \mbox{Neumann}.
\end{array} \right.$$
\end{prop}

\begin{proof} It  follows by Proposition \ref{CONORMAL} and a standard Tauberian argument
that there exists a two-term asymptotic expansion for some
principal coefficient. We may determine it from an integrated
version of the asymptotics where we integrate the left side
against a smooth function $\rho$ on $\partial \Omega$.  The
integrated version follows from (\ref{DSING}) in the Dirichlet
case and (\ref{NSING}) in the Neumann case.

\begin{lem} For any smooth $\rho$ on $\partial \Omega$,

$$\sum_{j: \lambda_j \leq \lambda}
\int_Y \rho |u_j^b(q)|^2 dA = \left\{\begin{array}{ll}
\int_{\partial \Omega} \rho dA  \lambda^{n + 2} +
O(\lambda^{n+1}), & \mbox{Dirichlet}
\\ & \\
 \int_{\partial \Omega} \rho dA \lambda^n + O(\lambda^{n-1}), & \mbox{Neumann}.
\end{array} \right.$$
\end{lem}

\end{proof}

We now improve the result in the generic case:

\begin{theo} Suppose that the set of loops at $q$ has measure $0$
in $B^*_q\partial \Omega$. Then

$$\sum_{j: \lambda_j \leq \lambda}
|u_j^b(q)|^2 = \left\{\begin{array}{ll} C \lambda^{n + 2} +
o(\lambda^{n+1}) , & \mbox{Dirichlet}
\\ & \\
C \lambda^n + o(\lambda^{n-1}), & \mbox{Neumann}.
\end{array} \right.$$

It follows that $||u_j^b||_{L^{\infty}(\partial \Omega)} =
o(\lambda^{\frac{n+1}{2}})$ in the Dirichlet case, resp.
$o(\lambda^{\frac{n-1}{2}})$ in the Neumann case. These results
are sharp.

\end{theo}

\begin{proof} The proof follows the outline of  that in \cite{SZ}.  We write the left
side as the boundary trace $E^b_{[0,\lambda]}(q,q)$ of the
spectral projections, and  define the remainders in the local Weyl
laws by
$$E^b_{[0,\lambda]}(q,q)= \left\{\begin{array}{ll} C \lambda^{n + 2} +
R_D(\lambda, q) , & \mbox{Dirichlet}
\\ & \\
C \lambda^n + R_N(\lambda, q), & \mbox{Neumann}.
\end{array} \right.$$

We first show that  if the set of billiard loops at $q \in
\partial \Omega$ has measure zero in $B^*_q\partial \Omega$,
 then given $\varepsilon>0$, we can find a ball $B$
centered at $q$ and a $\Lambda<\infty$ so that for $\lambda \geq
\Lambda$,
\begin{equation}\label{W1}
 \left\{
\begin{array}{l} |R_D(\lambda,q)|\le \varepsilon\lambda^{n + 1}, \, \, q\in B,\\ \\
|R_N(\lambda,q)|\le \varepsilon\lambda^{n - 1}, \, \, q\in B.
\end{array} \right.
\end{equation}

To prove this, we need some more notation. The loop-length
function on $B^*\partial \Omega$ is the lower semi-continuous
funciton  defined by
\begin{equation}\label{M1}
L^*(q, \eta)=\left\{ \begin{array}{l}\inf \{n>0: \, \pi \circ
\beta^n(q, \eta) = q \}, \\ \\
\infty, \;\; \mbox{if no such $n$ exists}. \end{array} \right.
\end{equation}
The set of loop directions at $q \in \partial \Omega$ is defined
by:
\begin{equation} {\mathcal L}_q = \{\eta \in B^*_q \partial \Omega \, 1/L^*(q, \eta)\ne 0\}.  \end{equation}

We  construct two semiclassical pseudodifferential operators $b(q,
D), B(q, D) = I - b(x, D)$ on $L^2(\partial \Omega)$ with the
property that $B$ is microsupported in the set where $L^*(q, \eta)
>> T$. We then study $E_{[0,\lambda]}^b(q,q)$ by writing it in the form:
\begin{equation} \label{MPU} E_{[0,\lambda]}^b(q,q,) = [(B + b) E_{\lambda}^b (B + b)](q,q). \end{equation}

We choose  $b\in C^\infty (S^{n-1})$  so that
\begin{equation} \label{SMALL} \int_{S^{n-1}}b(\eta)d\sigma(\eta)\le 1/T^2, \end{equation} and
$$|L^*(q, \eta)|\le 1/T, \quad \text{on} \, \, {\mathcal N}\times \text{supp }
B,$$ where ${\mathcal N}$ is a neighborhood of $q$. We then put
$$B(\xi)=1-b(\xi).$$ By  construction,
\begin{equation} \label{SMOOTH} E^b(t)B^*(q,q),\;\; BE^b(t,q,q) \in C^{\infty}(0,
T). \end{equation}

A calculation using the conormal singularity at $t = 0$ of $B
E^b(t, q,q), b E^b(t, q, q)$ shows that,  for $\lambda\ge 1$,

\begin{equation}\label{W6}
\bigl| R_N(\lambda, q) \bigr|\le
CT^{-1}\lambda^{n-1}+C_T\lambda^{n-2},
\end{equation}
where $C_T$ depends on $T$ but $C$ does not.  This of course
yields \eqref{W1}. The argument is similar for the Dirichlet case.

\end{proof}

\subsection{Applications to eigenfunctions}

We now use that bounds on eigenfunctions (or for spectral
projections for intervals of shrinking width) can be obtained from
the jump in the remainder:
\begin{equation}\label{eig3} \sum_{j: \lambda_j = \lambda} |u_j^b(q)|
=
\sqrt{R(\lambda,q)-R(\lambda-0,q)}.
\end{equation}
To complete the proof, we observe that for fixed $q \in
\partial \Omega$ and any $\varepsilon>0$ then  one
can find a neighborhood $\mathcal{N}_\varepsilon(q)$ of $q$ and an
$\Lambda_\varepsilon(q)$ so that when $\lambda\ge
\Lambda_\varepsilon(xq)$ and $y\in \mathcal{N}_\varepsilon(q)$ we
have $|R(\lambda,y)|\le \varepsilon \lambda^{n-1}$.  This implies
that $|u_j^b(y)|\le \varepsilon \lambda^{(n-1)/2}$ if $y\in
\mathcal{N}_\varepsilon(q)$ and $\lambda\ge
\Lambda_\varepsilon(q)$. Since $M$ is compact and since the
$\mathcal{N}_{\varepsilon}(q)$  form open cover of $M$, we may
choose a finite subcover and extract the largest
$\Lambda_{\varepsilon}(q)$. For this $\Lambda_\varepsilon$, we get
$$|u_j^b(y)|\le \varepsilon\lambda^{(n-1)/2},
\quad \lambda\ge \Lambda_{\varepsilon},$$ The
$o(\lambda^{(n-1)/2})$ bound follows
 since $\Lambda_{\varepsilon}$ depends only on $\varepsilon$.

\subsection{\label{EXAMPLES} Examples}

For generic metrics and/or boundaries, there are no recurrent
points \cite{SZ}. Moreover,  convex analytic domains in $\R^n$
never have have recurrent points on the boundary. Suppose to the
contrary that there exists $x_0 \in
\partial \Omega$ such that a positive measure of  geodesic rays
starting at $x_0$ return to $x_0$ at the same time. By
analyticity, all rays starting at $x_0$ return to $x_0$, and they
must all return at the same time. In particular,  boundary
geodesics have to return to $x_0$.  So do
 creeping rays which make a small angle to the boundary, and
which accumulate at  boundary geodesics. But the creeping rays
which approach a boundary geodesic are strictly shorter than the
limiting boundary geodesic in the Euclidean metric, so the return
time could not be the same. Therefore, the $L^{\infty}$ estimate
on boundary traces of eigenfunctions or spectral projections is
never sharp for convex analytic Euclidean domains.

On the other hand, the $L^{\infty}$ bound is sharp for the
Euclidean half-circle or for any sector of a circle. Indeed, the
invariant Neumann eigenfunctions (under rotations) of the  disc
achieve the interior $L^{\infty}$ bound at the center. On a
half-circle or sector, they remain Neumann eigenfunctions and
their boundary traces achieve the maximal $L^{\infty}$ bound
above. In the case of Dirichlet boundary conditions, the bounds
are also saturated on such domains (by taking the boundary trace
of the imaginary part of the eigenfunction transforming by $e^{i
\theta}$ under rotation by angle $\theta$).

Above,  we have mainly considered simply connected domains, but
the same questions may be posed for multiply connected plane
domains (for instance). It seems likely that any analytic (or even
piecewise analytic) domain which achieves the maximum $L^{\infty}$
bound  must be simply connected. This would generalize the
topological result of \cite{SZ} that the  maximum $L^{\infty}$
bound can only be achieved for metrics on the sphere.

 The question thus arises to find the maximal growth rate of
 $L^{\infty}$
 norms of boundary traces of eigenfunctions of smooth Euclidean domains,
 (or of convex analytic Euclidean domains),
and which domains achieve the bounds? In fact, the answer depends
on whether one studies relative sup-norms
$||u_j^b||_{L^{\infty}(\partial \Omega)} / ||u_j^b||_{L^2(\partial
\Omega)}$ or absolute sup-norms $||u_j^b||_{L^{\infty}(\partial
\Omega)}.$ Moreover, the sharpness depends on whether we consider
individual eigenfunctions or spectral projections for shrinking
intervals.

We do not know the answer to the question, but offer some
speculations. By recent results of Smith-Sogge \cite{SS},
following  earlier results of Grieser (in his unpublished PhD
thesis and in  \cite{G}; (see also \cite{S})), it is known that
whispering gallery Neumann modes of two-dimensional  convex
Euclidean domains saturate the $L^p$ norms with $2 \leq p \leq 8$
because they live in thin $\lambda^{-2/3}$ layers around the
boundary and therefore are of size $\lambda^{1/3}$ there. The
analogous problem in higher dimensions is still open.

However, they are not extremal for relative sup norms because  are
spread out all over the boundary. For instance, in the case of the
disc, the relative sup norm $||u_j^b||_{L^{\infty}(\partial
\Omega)} / ||u_j^b||_{L^2(\partial \Omega)}$  is equal to one.  A
better guess is that boundary traces of modes (or quasimodes)
associated to stable elliptic orbits are extremals for the
relative sup norm $||u_j^b||_{L^{\infty}(\partial \Omega)} /
||u_j^b||_{L^2(\partial \Omega)}$. They live in $\lambda^{-1/2}$
tubes around the orbits and therefore the modes have size
$\lambda^{(n-1)/4}$ along the orbit in dimension $n$. The boundary
trace is therefore concentrated in small balls of radius
$\lambda^{-1/2}$  around the bounce points of the orbit. So the
$L^2$ norm of the boundary trace should be $\sim 1$.

\end{document}